\documentclass[11pt]{amsart}

\usepackage{amsmath, amsthm, 
amsfonts, amssymb}

\swapnumbers

\newtheorem{Theorem}[subsection]{Theorem}
\newtheorem{Lemma}[subsection]{Lemma}
\newtheorem{Proposition}[subsection]{Proposition}

\newtheorem*{Claim}{Claim}
\newtheorem*{Cass}{Casselman's criterion} 
\newtheorem*{Main}{Main Theorem} 

\theoremstyle{definition}

\theoremstyle{remark}
\newtheorem{Remark}[subsection]{Remark}

\numberwithin{equation}{subsection}

\begin{document}

\title[Square integrability 
on $p$-adic symmetric spaces]
{Square integrability of 
representations \\ 
on $p$-adic symmetric spaces} 

\author{Shin-ichi Kato} 
\address{Department of Mathematics, Graduate 
 School of Science, 
 Kyoto University, 
 Kyoto 606--8502 Japan}
 \email{skato@math.kyoto-u.ac.jp}
 
\author{Keiji Takano} 
 \address{Department of Arts and Science, 
 Akashi National College of Technology, 
 679-3 Nishioka, Uozumi-cho, Akashi 674-8501 Japan}
 \email{takano@akashi.ac.jp}


\begin{abstract}A symmetric space analogue of 
Casselman's criterion for 
square integrability of representations of a 
$p$-adic group is established. 
It is described in terms of 
exponents of Jacquet modules along 
parabolic subgroups associated to 
the symmetric space. 
\end{abstract}

\maketitle

\section*{Introduction} 

Let ${\bf G}$ be a connected reductive group defined over 
a non-archimedean local field $F$ 
equipped with an involutive $F$-automorphism $\sigma: {\bf G}\to{\bf G}$, 
and ${\bf H}$ the subgroup of all $\sigma$-fixed points 
of ${\bf G}$. The quotient space $G/H$ where 
$G={\bf G}(F)$ and $H={\bf H}(F)$ 
is called a $p$-adic symmetric space. 
We are interested in representations of 
$G$ which can be realized in the space of functions on $G/H$. 
Such representations are often said to be {\it $H$-distinguished}. 
In this paper, we are concerned especially with 
{\it discrete series for $G/H$}; roughly speaking, the representations 
which have realizations in the space of 
square integrable functions on $G/H$. 

Let $(\pi,V)$ be a finitely generated admissible 
representation 
of $G$ which carries a non-zero 
$H$-invariant linear form $\lambda\in (V^*)^H$. 
We consider the functions $\varphi_{\lambda,v}$ on $G/H$ 
for $v\in V$ given by 
$$
\varphi_{\lambda,v}(g)=\langle \lambda,\pi(g^{-1})v\rangle\quad(g\in G).  
$$
Such functions are called 
{\it $H$-matrix coefficients of $\pi$ defined by $\lambda$}. 
Any non-trivial realization of $\pi$ in the space of functions on $G/H$ is 
formed by these $H$-matrix coefficients for some 
non-zero $\lambda\in (V^*)^H$. 
Note that $H$-matrix coefficients are not the 
matrix coefficients in the usual sense, but are 
{\it generalized} matrix coefficients, 
since $H$-invariant linear forms are not smooth in general. 
In a previous work \cite{KT}, we have studied representations 
whose $H$-matrix coefficients have compact support 
modulo $Z_GH$. Here $Z_G$ denotes the center of $G$. 
We have called such representations 
$(H,\lambda)$-relatively cuspidal, and 
established a criterion for 
$(H, \lambda)$-relative cuspidality of $\pi$ 
by using Jacquet modules (\cite[6.2]{KT}). 
In the present work, 
we shall deal with a different class of representations. 
For simplicity, suppose that $\pi$ has a unitary 
central character. Then 
$\bigl|\varphi_{\lambda,v}(\cdot)\bigr|$ 
is regarded as a function on $G/Z_GH$. We say that $\pi$ is 
{\it $H$-square integrable with respect to $\lambda$} if 
$\bigl|\varphi_{\lambda,v}(\cdot)\bigr|$ is square integrable 
on $G/Z_GH$ for all $v\in V$, namely, if 
$$
\int_{G/Z_GH} |\varphi_{\lambda,v}(g)|^2dg<\infty
$$
for all $v\in V$. 
We shall establish a 
criterion for $H$-square integrability of $\pi$ 
in this paper. 

Before stating our main theorem, 
let us recall Casselman's criterion for the 
usual 
square integrability. We say that $\pi$ is 
{\it square integrable} if all the usual matrix 
coefficients (defined by smooth linear forms) 
are square integrable on $G/Z_G$. 
For each parabolic subgroup $P$ of $G$ with the 
$F$-split component $A_P$, let $(\pi_P, V_P)$ be the 
normalized Jacquet module of 
$\pi$ along $P$ and ${\mathcal Exp}_{A_P}(\pi_P)$ 
the set of all quasi-characters $\chi$ of $A_P$ having 
non-zero generalized eigenvectors in $V_P$. 
Let $A_P^-$ and $A_P^1$ denote the dominant part of $A_P$ and 
the $\mathcal O_F$-points of $A_P$ respectively. 

\begin{Cass} {\rm{(\cite[4.4.6]{C})}} \label{Cass}  
The representation $\pi$ is square integrable if and 
only if for every parabolic subgroup $P$, the condition 
$\bigl|\chi(a)\bigr|<1$ 
holds for all $\chi\in{\mathcal Exp}_{A_P}(\pi_P)$ and 
$a\in A_P^-\setminus Z_GA^1_P$. 
\end{Cass} 

Also our criterion for $H$-square integrability is stated in terms of 
exponents of Jacquet modules. However we use only 
those along {\it $\sigma$-split parabolic subgroups} (see \ref{1.5}). 
In our previous work \cite{KT} (and also in Lagier \cite{L}), 
a canonical mapping 
$r_P: (V^*)^H\to (V_P^*)^{M\cap H}$ was 
introduced for each $\sigma$-split parabolic subgroup 
$P=MU$ (see \ref{3.2}). 
Now, for a given $\lambda\in (V^*)^H$ and a 
$\sigma$-split parabolic subgroup $P$ 
with the $(\sigma, F)$-split component $S_P$  (see \ref{1.5}), we put 
$$
{\mathcal Exp}_{S_P}\bigl(\pi_P,r_P(\lambda)\bigr)
=\left\{\chi\in {\mathcal Exp}_{S_P}(\pi_P)\left|\, 
\begin{matrix} r_P(\lambda)\neq 0 \text{ on the generalized} \\  
\text{$\chi$-eigenspace in } V_P\end{matrix}\right. \right\}. 
$$
The main theorem of this paper is the following (Theorem \ref{4.7}): 
\begin{Main} The representation $\pi$ is $H$-square integrable with 
respect to $\lambda$ if and only if 
for every $\sigma$-split parabolic subgroup $P$, 
the condition $\bigl|\chi(s)\bigr|<1$ holds for all 
$\chi\in {\mathcal Exp}_{S_P}\bigl(\pi_P,r_P(\lambda)\bigr)$ and 
$s\in S_P^-\setminus Z_GS_P^1$. 
\end{Main}  
This is an analogue, and even a generalization, of Casselman's criterion. 
Actually, if one applies the above theorem to the 
symmetric space $G/H=(G_1\times G_1)/G_1$ where the 
involution is 
the permutation of two factors 
(referred to as the {\it group case}), then one recovers 
Casselman's criterion for the group $G_1$. 

Let us summarize the contents of this paper. 
In Section 1 we prepare notation and several definitions used throughout. 
In Section 2, we recall the analogue of 
Cartan decomposition for $p$-adic symmetric spaces 
given by Benoist-Oh \cite{BO} and Delorme-S\'{e}cherre \cite{DS}. 
After that, we give 
two ingredients for the proof of the main theorem; 
a disjointness assertion (Proposition \ref{2.3}) and some 
volume estimate (Proposition \ref{2.6}). 
Section 3 is essentially a recollection of \cite[\S5]{KT} and \cite[\S2]{L} on 
the asymptotic behavior of $H$-matrix coefficients described 
by the mapping $r_P$. 
Section 4 is devoted to the proof of the main theorem. 
We give simple examples of $H$-square integrable 
representations in Section 5. 

\section{Notation and definitions}  

\subsection{Basic notation} \label{1.1} 
%{1.1}

\indent 

 Let $F$ be a non-archimedean local field with the 
absolute value $|\cdot|_F$. The valuation ring of $F$ is denoted by 
${\mathcal O}_F$ and the order of the residue field of $F$ by $q_F$. 
Throughout this paper, we assume that 
the residual characteristic of $F$ is not equal to $2$. 
Let $\nu_F: F^{\times}\to{\Bbb Z}$ denote the 
additive valuation defined by 
$$
|x|_F=q_F^{-\nu_F(x)}\quad (x\in F^{\times}). 
$$

Let $\mathbf{G}$ be a connected 
reductive group 
defined over $F$ and 
$\sigma$ an $F$-involution on $\mathbf{G}$. 
The $F$-subgroup $\{h\in\mathbf{G}\,|\,\sigma(h)=h\}$ consisting 
of all $\sigma$-fixed points of 
$\mathbf{G}$ is denoted by $\mathbf{H}$. Let 
$\mathbf{Z}$ be the $F$-split component of 
$\mathbf{G}$, that is, the largest 
$F$-split torus lying in the center of 
$\mathbf{G}$. Note that $\mathbf{Z}$ is 
$\sigma$-stable. 
We put 
$$
{\bf Z}_0=\left(\{z\in {\bf Z}\bigm| \sigma(z)=z^{-1}\}\right)^0,
$$
which we shall call 
{\it the $(\sigma,F)$-split component of ${\bf G}$}. Here and henceforth, 
$(\,\,)^0$ stands for the identity component in the 
Zariski topology. 

The group $\mathbf{G}(F)$ consisting of all the $F$-rational points 
of $\mathbf{G}$ is denoted by $G$.
Similarly, for any $F$-subgroup 
$\mathbf{R}$ of
$\mathbf{G}$, we shall write $R=\mathbf{R}(F)$. 

For a connected $F$-group ${\bf M}$, let $X^*({\bf M})$ (resp. 
$X^*({\bf M})_F$) denote the free $\mathbb Z$-module of 
rational (resp. $F$-rational) characters of ${\bf M}$. 
If ${\bf A}$ is an $F$-split torus, one has $X^*({\bf A})=X^*({\bf A})_F$. 
Let ${\bf M}$ be a connected reductive $F$-group with its 
$F$-split component ${\bf A}$. We put 
$$
{\frak a}_{\bf M}={\mathrm{Hom}}\left(X^*({\bf M})_F, \mathbb R\right). 
$$ 
The natural homomorphism $X^*({\bf M})_F\to X^*({\bf A})$ defined by 
restriction induces an isomorphism 
\begin{equation} \label{1.1.1}
{\mathrm{Hom}}\left(X^*({\bf A}), \mathbb R\right)\simeq 
{\frak a}_{\bf M}. 
\end{equation}
We define a homomorphism 
$$
\nu_M: M={\bf M}(F) \to {\frak a}_M
$$
by 
$$
\langle\alpha, \nu_M(m)\rangle =
\nu_F(m^{\alpha})
$$
for all $m\in M$ and 
$\alpha\in X^*({\bf M})_F$. 
We can define $\langle\alpha, \nu_M(m)\rangle $ also for 
$\alpha\in X^*({\bf A})$ through the 
identification \eqref{1.1.1}. The kernel of $\nu_M$ is denoted 
by $M^1$. Note that $A^1={\bf A}(\mathcal O_F)$ for 
an $F$-split torus ${\bf A}$. 

\subsection{$H$-matrix coefficients of representations} \label{1.2} 
%{1.2} 

\indent 

Let $C^{\infty}(G/H)$ denote the space of all 
smooth $\Bbb C$-valued 
functions on $G/H$ on which $G$ acts by 
left translation. 
A smooth representation $(\pi,V)$ of $G$ is said to be 
{\it $H$-distinguished} if $(V^*)^H\neq \{0\}$. 
Take a non-zero $\lambda\in(V^*)^H$ and 
consider the functions on $G$ given by 
$$
\varphi_{\lambda, v}(g)=\langle \lambda, \pi(g^{-1})v\rangle\quad
(g\in G)
$$
for $v\in V$. We call these functions 
{\it $H$-matrix coefficients of $\pi$ defined by $\lambda$}. 
Let us identify right $H$-invariant functions 
on $G$ with functions on $G/H$. Then, $H$-matrix coefficients 
belong to $C^{\infty}(G/H)$ 
and the mapping 
$$
T_{\lambda}: V\to C^{\infty}(G/H),\quad T_{\lambda}(v)=
\varphi_{\lambda,v}
$$
gives a $G$-morphism. Any realization of 
$V$ in $C^{\infty}(G/H)$ is determined by 
an $H$-invariant linear form in this way. 

Let $\omega_0$ be a 
quasi-character of $Z_0={\bf Z}_0(F)$. 
A smooth representation 
$(\pi,V)$ is called an $\omega_0$-representation if 
$Z_0$ acts on $V$ by the the character $\omega_0$. 
If $(\pi,V)$ is an $H$-distinguished $\omega_0$-representation, 
then for any $\lambda\in (V^*)^H$, the image 
$T_{\lambda}(V)$ of $V$ is contained in the space 
$C^{\infty}_{\omega_0}(G/H)$ consisting 
of functions $\varphi\in C^{\infty}(G/H)$ which satisfy 
$$
\varphi(zgH)=\omega_0(z)^{-1}\varphi(gH)\quad(z\in Z_0,\,gH\in G/H). 
$$

\subsection{$H$-square integrability} \label{1.3} 
%{1.3} 

\indent 

Both $G$ and $H$ are unimodular groups. 
Thus the quotient space 
$G/Z_0H$ carries a left $G$-invariant measure, which is 
denoted by 
$\int_{G/Z_0H}\,\, dg$. 
Let $\omega_0$  be a unitary character of $Z_0$. 
We define $L^2_{\omega_0}(G/H)$ to be the space of 
($L^2$-classes of) all functions 
$\varphi\in C_{\omega_0}^{\infty}(G/H)$ satisfying 
$$
\int_{G/Z_0H} |\varphi(g)|^2dg<\infty. 
$$
Take a non-zero $\lambda\in\bigl(V^*\bigr)^H$. 
We say that an $H$-distinguished 
$\omega_0$-representation 
$(\pi,V)$ is {\it $H$-square integrable with respect to 
$\lambda$} if the $H$-matrix coefficients $\varphi_{\lambda,v}$ defined by 
$\lambda$ are square integrable on $G/Z_0H$ for all $v\in V$, or 
equivalently, if 
$T_{\lambda}(V)$ is contained in 
$L^2_{\omega_0}(G/H)$. 
Note that this definition agrees with the one given in the 
Introduction, since 
$Z_GH/Z_0H$ is compact. Irreducible $H$-square integrable representations 
are said to be {\it in the discrete series for $G/H$}. 

In \cite{KT}, we have put the following 
definition: 
An $H$-distinguished admissible $\omega_0$-representation 
$(\pi,V)$ of $G$ is said to be {\it $(H,\lambda)$-relatively cuspidal} if 
all the $H$-matrix coefficients of $\pi$ defined by 
$\lambda$ are compactly supported modulo $Z_0H$. 
We gave examples of such representations in \cite[\S8]{KT}. 
It is clear from the definition that 
$(H,\lambda)$-relatively cuspidal $\omega_0$-representations 
are $H$-square integrable with respect to $\lambda$ provided that 
$\omega_0$ is unitary. 

\subsection{Tori and roots associated to symmetric spaces} \label{1.4} 

\indent 

We shall fix notation on tori, roots, and parabolic subgroups associated to 
the involution $\sigma$. For the reference, see \cite{HH} and also 
\cite[\S2]{KT}. 

A torus ${\bf S}$ is said to be $(\sigma,F)$-split if 
it is $F$-split and $\sigma(s)=s^{-1}$ for all $s\in{\bf S}$. 
We fix a maximal $(\sigma, F)$-split torus ${\bf S}_0$ of 
${\bf G}$ and a maximal $F$-split torus ${\bf A}_{\emptyset}$ containing 
${\bf S}_0$. Then ${\bf A}_{\emptyset}$ is 
necessarily $\sigma$-stable, so $\sigma$ acts 
naturally on $X^*({\bf A}_{\emptyset})$. 

Let $\Phi\subset {\rm X}^*({\bf A}_{\emptyset})$ be the root system of 
$({\bf G}, {\bf A}_{\emptyset})$. 
It is $\sigma$-stable. We choose a $\sigma$-basis 
$\Delta$ of $\Phi$ that has the property 
$$
\alpha>0,\,\sigma(\alpha)\neq 
\alpha\,\Longrightarrow\,\sigma({\alpha})<0
$$
in the corresponding order. The subset of all $\sigma$-fixed roots in 
${\Phi}$ (resp. $\Delta$) is 
denoted by ${\Phi}_{\sigma}$ (resp. $\Delta_{\sigma}$). 

Let 
$p: {\rm X}^*({\bf A}_{\emptyset})\to{\rm X}^*({\bf S}_0)$ be the 
homomorphism defined by restriction to ${\bf S}_0$. It is surjective 
and its kernel coincides with 
the submodule of all $\sigma$-fixed elements 
of ${\rm X}^*({\bf A}_{\emptyset})$. Let us put 
$$
\overline{\Phi}=p(\Phi)\setminus\{0\}\,\,\Bigl(=p({\Phi}\setminus {\Phi}_{\sigma})\Bigr). 
$$ 
It is well-known that $\overline{\Phi}$ is a root system in 
${\rm X}^*({\bf S}_0)$ with a basis 
$$
\overline{\Delta}=
p(\Delta)\setminus \{0\}\,\,
\Bigl( =p({\Delta}\setminus {\Delta}_{\sigma})\Bigr). 
$$
For each subset $\overline{I}$ of $\overline{\Delta}$, 
we consider the subset 
$$ 
[\overline{I}]:=
\left(p^{-1}\left(\overline{I}\right)\cap\Delta\right)\cup\Delta_{\sigma}
$$ 
of $\Delta$. In this paper we say that 
a subset of $\Delta$ is {\it $\sigma$-split} if it is 
of the form $[\overline{I}]$ for some $\overline{I}\subset
\overline{\Delta}$. This terminology agrees with 
that in \cite[2.3]{KT}. 
The correspondence 
$\overline{I}\mapsto [\overline{I}]$ is an inclusion-preserving 
bijection between subsets of $\overline{\Delta}$ and 
$\sigma$-split subsets of $\Delta$. 
The inverse of this correspondence 
is given by 
$I\mapsto p({I}\setminus {\Delta}_{\sigma})$ 
for a $\sigma$-split subset $I$ of $\Delta$. 

Note that maximal proper $\sigma$-split subsets of 
$\Delta$ are written in the form 
$[\overline{\Delta}\setminus\{\overline{\alpha}\}]$ 
for some $\overline{\alpha}\in \overline{\Delta}$. 

\subsection{Parabolic subgroups associated to 
symmetric spaces} \label{1.5}%{1.5} 

\indent 

A parabolic $F$-subgroup ${\bf P}$ of ${\bf G}$ is said to be 
$\sigma$-split if ${\bf P}$ and $\sigma({\bf P})$ are opposite. 
In such a case, we always take ${\bf M}={\bf P}\cap\sigma({\bf P})$ 
for a ($\sigma$-stable) Levi subgroup of ${\bf P}$. 

Let ${\bf P}_{\emptyset}$ be the minimal parabolic $F$-subgroup of 
${\bf G}$ corresponding to the choice of ${\Delta}$ as in \ref{1.4}. 
The centralizer $Z_{\bf G}({\bf A}_{\emptyset})$ of 
${\bf A}_{\emptyset}$ in ${\bf G}$ is denoted by ${\bf M}_{\emptyset}$, 
which is a Levi subgroup of ${\bf P}_{\emptyset}$. 
Put ${\bf M}_0=Z_{\bf G}({\bf S}_0)$ and 
${\bf P}_0={\bf P}_{\emptyset}{\bf M}_0$. 
Then ${\bf P}_0$ is a minimal $\sigma$-split parabolic 
$F$-subgroup of ${\bf G}$ with a $\sigma$-stable Levi subgroup 
${\bf M}_0$. 
Let ${\bf U}_0$ be the unipotent radical of ${\bf P}_0$. 

For each subset $I$ of $\Delta$, let ${\bf P}_{{I}}$ denote 
the standard parabolic $F$-subgroup 
(that contains ${\bf P}_{\emptyset}$) 
of ${\bf G}$ corresponding to ${I}$. 
If $I$ is a $\sigma$-split subset of $\Delta$, 
then ${\bf P}_{{I}}$ is $\sigma$-split. 
Standard $\sigma$-split parabolic $F$-subgroups of ${\bf G}$ 
correspond to $\sigma$-split subsets of $\Delta$ in this way. 
Furthermore, it is remarked in \cite[2.5]{KT} that 
any $\sigma$-split parabolic 
$F$-subgroup 
of ${\bf G}$ is of the form $\gamma^{-1}{\bf P}_{{I}}\gamma$ 
for some 
$\sigma$-split $I\subset\Delta$ and $\gamma\in({\bf M}_0{\bf H})(F)$. 
For each $\sigma$-split $I\subset\Delta$, we put 
${\bf M}_{{I}}={\bf P}_{{I}}\cap \sigma({\bf P}_{{I}})$. 
Note that ${\bf P}_{\Delta_{\sigma}}$ and 
${\bf M}_{\Delta_{\sigma}}$ 
coincide with ${\bf P}_0$ and ${\bf M}_0$ respectively. 
Let ${\bf U}_{{I}}$ be the unipotent radical of ${\bf P}_{{I}}$, 
so that we have a Levi decomposition 
${\bf P}_{{I}}={\bf M}_{{I}}{\bf U}_{{I}}$. 
The $F$-split component (resp. $(\sigma,F)$-split component) of 
${\bf M}_{{I}}$ is denoted by ${\bf A}_I$ (resp. ${\bf S}_I$). 
These are also called the $F$-split and $(\sigma,F)$-split 
component of ${\bf P}_I$ respectively. 

We can describe $\sigma$-split parabolic subgroups 
and their $(\sigma,F)$-split components also by using 
subsets of the restricted basis $\overline{\Delta}$ 
through the bijective correspondence 
$\overline{I}\leftrightarrow I= [\overline{I}]$ in \ref{1.4}. 
We shall occasionally use notation based on 
restricted roots if it is convenient. 
For each subset $\overline{I}$ of $\overline{\Delta}$, we put 
$$
{\bf S}_{\overline{I}}=
\left(\bigcap_{\overline{\alpha}\in \overline{I}} 
\ker\left(\overline{\alpha}: {\bf S}_0\to{\bf G}_{\rm m}\right)
\right)^0. 
$$
It is easy to see the following equalities: 
\begin{equation} \label{1.5.1} 
{\bf S}_0={\bf S}_{\overline{I}}\cdot{\bf S}_{\overline{\Delta}\setminus 
\overline{I}},\quad 
{\bf S}_{\overline{I}}\cap{\bf S}_{\overline{\Delta}\setminus \overline{I}}=
{\bf S}_{\overline{\Delta}}={\bf Z}_0. 
\end{equation}  
Observe that if $\overline{I}\subset\overline{\Delta}$ 
corresponds to a $\sigma$-split subset $I\subset\Delta$, then 
${\bf S}_{\overline{I}}$ coincides with ${\bf S}_I$. 

For a positive real number $\varepsilon\leqq 1$ and a 
$\sigma$-split subset $I$ of $\Delta$, 
we put 
$$
S_I^-(\varepsilon)=
\left\{ s\in S_I\bigm| |s^{{\alpha}}|_F\leqq\varepsilon\,\,
({\alpha}\in
{\Delta}\setminus {I})\right\},  
$$ 
$$
S_{0, I}^-(\varepsilon)=
\left\{ s\in S_0\left|\, 
\begin{matrix}  |s^{{\alpha}}|_F\leqq\varepsilon\,\,({\alpha}\in
{\Delta}\setminus {I}),\\ 
|s^{{\alpha}}|_F\leqq 1 \,\,({\alpha}\in {I})\end{matrix}\right.  
\right\},  
$$
and 
$$
{}_{I}S_0^-(\varepsilon)=
\left\{ s\in S_0\left|\, 
\begin{matrix} 
|s^{{\alpha}}|_F\leqq\varepsilon\,\,({\alpha}\in
{\Delta}\setminus {I}),\\ 
\varepsilon<|s^{{\alpha}}|_F\leqq 1\,({\alpha}\in 
{I})\end{matrix}\right. 
\right\}. 
$$
For example, $S_I^-(\varepsilon)$ can be written also as 
\begin{equation} \label{1.5.2} 
S_I^-(\varepsilon)=
\left\{ s\in S_I=S_{\overline{I}}\,\Bigl| \bigl|s^{{\overline{\alpha}}}\bigr|_F\leqq\varepsilon\,\,
\left(\overline{\alpha}\in
\overline{\Delta}\setminus \overline{I}\right)\Bigr.\right\}
\end{equation} 
if $I$ corresponds to $\overline{I}\subset\overline{\Delta}$ as in \ref{1.4}. 
We abbreviate $S_I^-(1)$ as $S_I^-$ and write 
$$
S_0^-\,\left(=S_{\Delta_{\sigma}}^-(1)\right)
\,=\left\{ s\in S_0\Bigm| 
 |s^{{\alpha}}|_F\leqq 1\,({\alpha}\in
{\Delta})\,
\right\}.
$$ 
It is obvious from the definition that 
$$
S_I^-(\varepsilon)\subset 
{}_IS_0^-(\varepsilon)\subset 
S_{0, I}^-(\varepsilon)
\subset 
S_0^-, 
$$
and that 
$$
\varepsilon\leqq \varepsilon'\,\Longrightarrow 
S_{0, I}^-(\varepsilon)\subset 
S_{0, I}^-(\varepsilon'). 
$$

\begin{Lemma}\label{1.6} %{1.6} 
Let $I$ be a $\sigma$-split subset of $\Delta$. 

$(1)$ For any positive real number 
$\varepsilon\leqq 1$, 
there exists a positive real number 
$\varepsilon'\leqq 1$ 
such that 
$$
S_{0,I}^-(\varepsilon')\subset S_I^-(\varepsilon)\cdot S_0^-. 
$$

$(2)$ For any positive real number 
$\varepsilon\leqq 1$, 
there exist a positive real number 
$\varepsilon'$ and finitely many elements $t_1,\cdots, t_k$ of $S_0^-$ 
such that 
$$
{}_I S_0^-(\varepsilon)\subset 
\bigcup_{i} S_I^-(\varepsilon') t_iZ_0S_0^1. 
$$
Moreover one can take $\varepsilon'\leqq 1$ if $\varepsilon$ is 
sufficiently small. 

$(3)$ For any positive real number $\varepsilon\leqq 1$, 
one has a decomposition 
$$
S_0^-=\bigcup_{
I\subset\Delta:
\sigma{\text{-}\rm{split}}}
{}_IS_0^-(\varepsilon) \,\,\text{$($disjoint$)$} 
$$
where $I$ ranges over all $\sigma$-split subsets of $\Delta$. 
\end{Lemma} 

\begin{proof}  
The proof of (1) is exactly the same as that of \cite[4.3.1]{C}. 
Regard the union in (3) as 
$$
\bigcup_{\overline{I}\subset\overline{\Delta}}
 {}_{I}S_0^-(\varepsilon)\quad (I=[\overline{I}]) 
$$
where $\overline{I}$ ranges over all subsets of $\overline{\Delta}$ 
including the empty set. Then the assertion of (3) can be seen 
by the same way as in \cite[remark preceding 4.3.4]{C}. 
For (2), let $\overline{I}$ be the subset of $\overline{\Delta}$ corresponding 
to $I$. First observe that 
$S_{\overline{I}}\cdot S_{\overline{\Delta}\setminus
\overline{I}}$ is 
of finite index in $S_0$ by \eqref{1.5.1}. We can 
take a finite set $\Gamma_{I}$ of representatives of 
$S_0/(S_{\overline{I}}\cdot S_{\overline{\Delta}\setminus
\overline{I}})$ from $S_0^-$. We put 
$$
c=c_{{I}}=\min_{\gamma\in\Gamma_{{I}}, 
{\overline{\alpha}}\in \overline{\Delta}\setminus \overline{I}} 
\bigl( \left|\gamma^{\overline{\alpha}}\right|_F\bigr) 
$$
and take $t_1,\dots, t_k\in S_0^-$ so that $\bigcup_i t_iZ_0S_0^1$ 
contains the subset 
$$
\left\{ s\in S_0^-\left|\, \begin{matrix} 
\varepsilon<\left| s^{\overline{\alpha}}
\right|_F\leqq 1\,\,(
\overline{\alpha}\in \overline{I}),\\
c\leqq \left| s^{\overline{\alpha}}\right|_F\leqq 1\,\,
(\overline{\alpha}\in\overline{\Delta}\setminus \overline{I})
\end{matrix}\right. \right\}. 
$$
Now, let us write $s\in {}_I S_0^-(\varepsilon)$ as 
$s=s_1s_2\gamma$ with $s_1\in S_I=S_{\overline{I}}$, $s_2\in 
S_{\overline{\Delta}\setminus
\overline{I}}$, and 
$\gamma\in\Gamma_{{I}}$. We show that 
$s_1\in S_I^-(\varepsilon')$ for some $\varepsilon'$ and that 
$s_2\gamma\in \bigcup_i t_iZ_0S_0^1$. 
For any $\overline{\alpha}\in\overline{\Delta}\setminus 
\overline{I}$, we have 
$$
\left| s_1^{\overline{\alpha}}\right|_F\cdot c\leqq 
\left| (s_1\gamma)^{\overline{\alpha}}\right|_F
=\left| s^{\overline{\alpha}}\right|_F\leqq \varepsilon. 
$$ 
Therefore $s_1$ belongs to $S_I^-(\varepsilon')$ for 
$\varepsilon'=\varepsilon c^{-1}$ (see \eqref{1.5.2}). 
Note that $\varepsilon'\leqq 1$ 
if $\varepsilon\leqq c=c_{{I}}$. On the other hand, 
we have 
$$
\varepsilon<\left| s^{\overline{\alpha}}\right|_F=
\left| (s_2\gamma)^{\overline{\alpha}}\right|_F
\leqq 1 
$$
for each $\overline{\alpha}\in \overline{I}$, 
while 
$$
c\leqq\left| \gamma^{\overline{\alpha}}\right|_F=
\left| (s_2\gamma)^{\overline{\alpha}}\right|_F\leqq 1
$$
for each $\overline{\alpha}\in\overline{\Delta}\setminus \overline{I}$. 
This completes the proof. 
\end{proof}  

\section{Relative Cartan decomposition}

In this section we recall the analogue of 
Cartan decomposition for $p$-adic symmetric spaces 
given by Benoist-Oh \cite{BO} and 
Delorme-S\'{e}cherre \cite{DS}. 
Concerning this decomposition, 
we shall give a disjointness 
result (Proposition \ref{2.3}) and some 
volume estimate (Proposition \ref{2.6}). 
These will be key ingredients for 
the proof of the main theorem in Section 4. 

\subsection{The analogue of Cartan decomposition}\label{2.1} 

\indent 

Let us fix the data $({\bf S}_0, {\bf A}_{\emptyset}, {\Delta})$ 
as in \ref{1.4}. We shall write $S_I^+(\varepsilon)$, 
$S_{0, I}^+(\varepsilon)$, and 
${}_{I}S_0^+(\varepsilon)$ 
for 
the set of all elements $s$ such that $s^{-1}$
belong to $S_I^-(\varepsilon)$, 
$S_{0, I}^-(\varepsilon)$, and 
${}_{I}S_0^-(\varepsilon)$ respectively. 
We recall the analogue of 
Cartan decomposition 
given in 
\cite{BO} and \cite{DS}. Let us state it in 
the following form as in 
\cite[3.4]{KT}: 
There exists a compact subset $\Omega$ of $G$ and 
a finite subset 
$\Gamma$ of $\left(\mathbf{M}_0\mathbf{H}\right)(F)$ 
such that 
$$
G=\Omega S_{0}^+\Gamma H. 
$$

\subsection{Groups with Iwahori factorization} \label{2.2} 

\indent 

Let $K$ be a $\sigma$-stable open compact subgroup of $G$ and 
${\bf P}={\bf M}{\bf U}$ a 
$\sigma$-split parabolic subgroup of ${\bf G}$ 
standard with respect to ${\Delta}$. 
We put ${\bf P}^-=\sigma({\bf P})$ and ${\bf U}^-=\sigma({\bf U})$, 
so that ${\bf P}\cap{\bf P}^-={\bf M}$ and 
${\bf P}^-={\bf M}{\bf U^-}$. We say that $K$ has the 
{\it Iwahori factorization with respect to 
${\bf P}$} if the product map 
induces a bijection 
$$
U^-_K\times M_K\times U_K\simeq K, 
$$
where $U^-_K=U^-\cap K$, $M_K=M\cap K$, and $U_K=U\cap K$. 

Let $K_{\rm max}$ be an ${\bf A}_{\emptyset}$-good maximal 
compact subgroup of $G$. We take a $\sigma$-stable 
open compact subgroup $K_0$, 
having the Iwahori factorization with respect to ${\bf P}_0$, 
such that $S_0^1\subset K_0\subset K_{\rm max}$. 
For example, 
it suffices to take $K_0=B\cap \sigma(B)$ for an Iwahori 
subgroup $B$ contained in $K_{\rm max}$. 
Choose a finite set $\Xi\subset G$ 
such that $\Omega\subset \bigcup_{\xi\in\Xi} \xi K_0$. 
We may decompose $G$ as 
$$
G=\bigcup_{\xi\in\Xi,\gamma\in\Gamma}\bigcup_{\dot{s}\in S_0^+/S_0^1} 
\xi K_0 \dot{s}\gamma H, 
$$
hence 
\begin{equation}\label{2.2.1}
G/Z_0H
=
\bigcup_{\xi,\gamma}\bigcup_{\dot{s}\in S_{0}^+/Z_0S_0^1} 
\xi K_0 \dot{s}\gamma H/H. 
\end{equation} 
It is not known whether these are disjoint. However 
we can assert the following.  
\begin{Proposition}\label{2.3} 
For each $\gamma\in({\bf M}_{0}{\bf H})(F)$, 
the union 
$$
\bigcup_{\dot{s}\in S_0^+/S_0^1} 
K_0 \dot{s}\gamma H
$$
is disjoint. 
\end{Proposition} 
\begin{proof}  We use the mapping $\tau: G\to G$ given by 
$\tau(g)=g\sigma(g)^{-1}$ 
and the action 
$$
(g, x)\mapsto g*x=gx\sigma(g)^{-1} \quad(g\in G,\,x\in\tau(G)) 
$$
of $G$ on $\tau(G)$ to study cosets modulo $H$. 
Suppose that $K_0 s_1\gamma H=K_0 s_2\gamma H$ for 
$s_1$, $s_2\in S_{0}^+$. Applying $\tau$, we have 
$$
K_0*(s_1^2m_{\gamma})=K_0*(s_2^2m_{\gamma}) 
$$
where $m_{\gamma}=\tau(\gamma)\in M_{0}$. 
There is an element $k\in K_0$ such that $k*(s_1^2m_{\gamma})=
s_2^2m_{\gamma}$. Thus 
$$
k s_1^2m_{\gamma}=s_2^2m_{\gamma}\sigma(k). 
$$
Using the Iwahori factorization, 
write $k$ and $\sigma(k)\in K_0$ as 
$$
k=u_1^-m_1u_1,\quad \sigma(k)=u_2^-m_2u_2
$$
where $u_i^-\in 
(U_0^-)_{K_0}$, $m_i\in 
(M_0)_{K_0}$, and $\,u_i\in 
(U_0)_{K_0}$. 
Then we have 
$$
u_1^-m_1u_1s_1^2m_{\gamma}=s_2^2m_{\gamma}u_2^-m_2u_2, 
$$
hence 
\begin{align} 
u^-_1&\cdot (m_1s_1^2m_{\gamma})\cdot 
(m_{\gamma}^{-1}(s_1^2)^{-1}u_1s_1^2m_{\gamma})\notag \\ 
&=(m_{\gamma}s_2^2u_2^-(s_2^2)^{-1}m_{\gamma}^{-1})\cdot 
(m_{\gamma}s_2^2m_2)\cdot u_2. \notag
\end{align} 
By the uniqueness of expressions in $U_0^-M_0U_0$, we must have 
\begin{equation} \label{2.3.1} 
m_1s_1^2m_{\gamma}=m_{\gamma}s_2^2m_2. 
\end{equation} 
Now we use the usual Cartan decomposition (see e.g. \cite[1.1(4)]{W}) for 
$M_0$ to write 
$m_{\gamma}\in M_0$ as 
$$
m_{\gamma}=k_1ak_2\quad \bigl(k_1,k_2\in M_0\cap K_{\mathrm{max}},\,
a\in \overline{M}_{\emptyset}^{\Delta_{\sigma}, +}\bigr). 
$$
Here $\overline{M}_{\emptyset}^{\Delta_{\sigma}, +}$ denotes 
the set of elements $m$ of $M_{\emptyset}$ satisfying 
$$
\langle\alpha,\,\nu_{M_{\emptyset}}(m)\rangle
\leqq 0
$$ 
for all $\alpha\in\Delta_{\sigma}$. Note that 
$m_1$, $m_2\in M_0\cap K_0\subset M_0\cap K_{\mathrm{max}}$, and 
$S_0^+\subset \overline{M}_{\emptyset}^{\Delta_{\sigma}, +}$. 
Since $s_1^2$ and $s_2^2$ are central in $M_0$, we have 
$$
m_1k_1\cdot (s_1^2a)\cdot k_2=k_1\cdot (s_2^2a)\cdot 
k_2m_2 
$$
by \eqref{2.3.1}. From the uniqueness of Cartan decomposition, 
we may conclude that 
$$
s_1^2a\equiv s_2^2a\mod M_{\emptyset}^1. 
$$
Since $M_{\emptyset}^1\cap S_0=S_0^1$, 
we have 
$s_1^2\equiv s_2^2\mod S_0^1$, 
and in turn, $s_1\equiv s_2\mod S_0^1$. 
\end{proof} 

\subsection{Some volume computation}\label{2.4} 

\indent 

Let us fix a left $G$-invariant measure 
$\int_{G/H}\,\, dg$ on the quotient space $G/H$. 
For each open compact subgroup $K$ of $G$ and 
an element $a\in G$, the set $K aH/H$ is open and compact in $G/H$. Let 
$\int_{K aH/H}\,\,dg$ denote the restriction 
of $\int_{G/H}\,\,dg$ to the open subset $K aH/H$. 

We have to evaluate the volume of $K_0 s\gamma H/H$. 
For this purpose we first need to take a 
$\sigma$-stable open 
compact subgroup $K$ from {\it an adapted family} 
given in \cite[4.3]{KT}: Besides the Iwahori factorization, we need 
a good filtration (\cite[4.3 (2)]{KT}) inside 
$K$ (and its subgroups) to use the result 
\cite[4.6]{KT}. Then we can compute the volume of 
$KsH/H$ for such a $K$ and $s\in S_0^+$. 

Let $\delta_P$ denote the modulus character of 
$P={\bf P}(F)$ for a parabolic $F$-subgroup ${\bf P}$ of 
${\bf G}$. 

\begin{Lemma}\label{2.5} 
If $K$ is a $\sigma$-stable 
open compact subgroup 
as above, then 
$$
{\rm vol}(K s H/H)=\delta_{P_0}(s)\cdot{\rm vol}(KeH/H) 
$$ 
for all $s\in S_0^+$. 
\end{Lemma} 
\begin{proof}  
We use the mapping $\tau$ to identify 
$\tau(G)$ with $G/H$ and transport the left $G$-invariant measure 
on $G/H$ to the measure on $\tau(G)$ invariant under the $*$-action of $G$. 
The subset 
$KsH/H$ is identified with 
$(Ks)*e$ and 
$$
{\rm vol}\bigl((Ks)*e\bigr)={\rm vol}\bigl(s*((s^{-1}Ks)*e)\bigr)
={\rm vol}\bigl((s^{-1}Ks)*e\bigr)
$$ 
by the invariance under $*$. Put $K_s=s^{-1}Ks$. The 
Iwahori factorization 
$$
K=U_{0,K}^-M_{0,K}U_{0,K}
$$
of $K$ with respect to ${\bf P}_0$ implies 
$$
K_s=(s^{-1}U_{0,K}^-s)\cdot 
M_{0,K}\cdot (s^{-1}U_{0,K}s) 
$$
for $s\in S_0^+$ with 
$$
s^{-1}U_{0,K}^-s\supset 
U_{0,K}^-, \quad 
s^{-1}U_{0,K}s\subset U_{0,K}. 
$$
Put $K_s'=K_s\cap \sigma(K_s)\,\bigl(=
s^{-1}Ks\cap sKs^{-1}\bigr)$. It is $\sigma$-stable and has 
a factorization 
$$
K_s'=(sU_{0,K}^-s^{-1})M_{0,K}(s^{-1}U_{0,K}s).  
$$
We have $K'_s\subset K$. 
So we may apply \cite[4.6]{KT} to the group $K_s'$, 
which asserts that 
$$
s^{-1}U_{0,K}s\subset 
(sU_{0,K}^-s^{-1})M_{0,K}H. 
$$
As a result, we have 
\begin{align} 
K_s&*e=(s^{-1}U_{0,K}^-s)
M_{0,K}(s^{-1}U_{0,K}s)*e \notag\\ 
&\subset (s^{-1}U_{0,K}^-s)
M_{0,K}(sU_{0,K}^-s^{-1})M_{0,K}*e
=(s^{-1}U_{0,K}^-s)
M_{0,K}*e.\notag 
\end{align}
We also note that 
$(s^{-1}U_{0,K}^-s)
M_{0,K}=P_0^-\cap K_s$. 
This shows that the orbit 
$(P_0^-\cap K_s)*e$ of the subgroup 
$P_0^-\cap K_s$ of $K_s$ coincides with 
the whole $K_s$-orbit 
$K_s*e$. Now, 
fix Haar measures $du^-$ and $dm$ on $s^{-1}U_{0,K}^-s$ and 
$M_{0,K}$ respectively. The $(P_0^-\cap K_s)$-invariant measure 
$$
f\mapsto \int_{s^{-1}U_{0,K}^-s}\int_{M_{0,K}}f\bigl(u^-*(m*e)
\bigr)dmdu^-
$$
on the orbit $K_s*e$ 
has to be $K_s$-invariant. Consequently, 
the volume of $K_s*e$ is proportional to 
$$
{\rm vol}(s^{-1}U_{0,K}^-s)\cdot
{\rm vol}(M_{0,K})
={\rm vol}(U^-_{0,K})\cdot\delta_{P_0}(s)\cdot
{\rm vol}(M_{0,K}),   
$$
which is a constant multiple of $\delta_{P_0}(s)$. 
The constant turns out to be 
the volume of $KeH/H$, if we consider $s=e$. 
\end{proof}  

\begin{Proposition} \label{2.6} 
Let $K$ be an arbitrary open compact subgroup of $G$. 
For each $\gamma\in \Gamma$, 
there exist positive real constants $c_1$ and 
$c_2$ 
such that 
$$
c_1\cdot\delta_{P_0}(s)
\leqq {\rm vol}(K s\gamma H/H)
\leqq c_2\cdot\delta_{P_0}(s)
$$
for all $s\in S_0^+$. 
\end{Proposition} 
\begin{proof}  
We put 
$K'=\gamma^{-1} K\gamma$, 
${\bf S}_0'=\gamma^{-1} {\bf S}_0\gamma$, 
and ${\bf P}_0'=\gamma^{-1} {\bf P}_0\gamma$. 
Then ${\bf P}_0'$ is a $\sigma$-split 
parabolic subgroup 
with the $(\sigma,F)$-split component ${\bf S}_0'$. 
We have 
$$
{\rm vol}(K s\gamma H/H)
={\rm vol}(\gamma K' \gamma^{-1}s\gamma H/H)
={\rm vol}(K's'H/H) 
$$ 
where $s'=\gamma^{-1}s\gamma\in (S_0')^+
=\gamma^{-1}S_0^+\gamma$. Take a member $K''$ from the 
adapted family (corresponding to the $\gamma$-conjugated 
data) such that $K''\subset K'$. 
Then we have 
$$
{\rm vol}(K'' s' H/H)\leqq 
{\rm vol}(K' s' H/H)\leqq 
[K':K'']\cdot{\rm vol}(K'' s' H/H). 
$$
Thus the claim follows from Lemma \ref{2.5}. 
\end{proof} 

\section{Asymptotic behavior of $H$-matrix coefficients} 

In this section we shall describe 
asymptotic behavior of $H$-matrix coefficients 
through the mapping $r_P$ (see \ref{3.2}). 
This section is 
essentially a recollection of \cite[\S5]{KT} and \cite[\S2]{L}. 

From now on, we shall 
briefly say that $P$ is a $\sigma$-split 
parabolic subgroup of $G$ if $P$ is the 
group of $F$-points of a $\sigma$-split 
parabolic $F$-subgroup $\mathbf{P}={\bf P}_{{I}}$ 
of $\mathbf{G}$ etc, by abuse of terminology. 

\subsection{Normalized Jacquet modules}%{3.1} 

\indent 

For an admissible representation 
$(\pi,V)$ of $G$ and a parabolic subgroup $P=MU$ of $G$, 
the normalized 
Jacquet module of $\pi$ along $P$ is denoted by $(\pi_P,V_P)$: 
The space 
$V_P$ is defined as the quotient 
$V/V(U)$ where $V(U)$ is the subspace of 
$V$ spanned by 
all the elements of the form $\pi(u)v-v$ 
($u\in U$, $v\in V$). The action $\pi_P$ of 
$M$ on $V_P$ is normalized so that
$$
\pi_P(m)j_P(v)=
\delta_P^{-1/2}(m)j_P(\pi(m)v)
$$
for $m\in M$ and $v\in V$ where 
$j_P:V\to V_P$ denotes 
the  canonical projection. 

\subsection{The mapping $r_P$} \label{3.2} 

\indent 

When $P=MU$ is a $\sigma$-split parabolic subgroup, 
we have defined in \cite{KT} a linear mapping 
$$
r_P: (V^*)^H\to (V_P^*)^{M\cap H}
$$
between the spaces of invariant linear forms. 
If $v\in V$ is a canonical lifting (in the sense of \cite[\S 4]{C}) of 
$\bar{v}\in V_P$ with 
respect to a suitable $\sigma$-stable open compact subgroup 
(a member in an {\it adapted family} in \cite[4.3]{KT}), 
then 
$r_P(\lambda)$ for $\lambda\in(V^*)^H$ is 
defined by
\begin{equation} \label{3.2.1} 
\langle r_P(\lambda), \bar{v}\rangle 
=\langle \lambda, v\rangle 
\end{equation} 
(see \cite[5.3~(2) and 5.4]{KT}). The same mapping was 
constructed independently by N. Lagier \cite{L} 
in a different manner. 
P. Delorme extended the construction of 
such mappings to any smooth representations in \cite{D}. 

In \cite[\S5]{KT}, we gave asymptotic behavior of 
$H$-matrix coefficients through the mapping $r_P$. 
The result \cite[Proposition 5.5]{KT} can be extended to the next proposition. 
This is a generalization of Casselman's result \cite[4.3.3]{C} to symmetric spaces, and 
was already proved essentially in \cite[Th\'eor\`eme 2]{L}. 
The proof provided by \cite{L} invokes Casselman's result itself. 
We shall give a proof which do not rely on \cite[loc.cit]{C} 
(but follow the lines similar to \cite{C} as we did in \cite[\S 5]{KT}). 
It would yield \cite[loc.cit]{C} when we apply this to the {\it group case}. 

\begin{Proposition}\label{3.3} 
Let $I$ be a $\sigma$-split subset of $\Delta$ and $P=P_{{I}}$ 
the corresponding $\sigma$-split parabolic subgroup with 
the $(\sigma, F)$-split component $S=S_I$. Let $(\pi,V)$ be   
an admissible representation of $G$ and 
$V_1\subset V$ a finite dimensional subspace. 
Then there exists 
a positive real number $\varepsilon
=\varepsilon_I\leqq 1$ such 
that  
$$
\langle\lambda,\pi(s)v\rangle
=
\delta_P^{1/2}(s)
\langle r_P(\lambda),\pi_P(s)j_P(v)\rangle
$$
for all $s\in S_{0, I}^-(\varepsilon)$, $v\in V_1$, and $\lambda\in(V^*)^H$. 
\end{Proposition} 
\begin{proof}  
For a compact subgroup $K$ of $G$, let $p_K: V\rightarrow V^K$ denote the 
projection defined by $p_K(v)=\int_K\pi(k)v\,dk$. For an open compact 
subgroup $K$ of $G$ (from the adapted family), 
there is a positive real number $\varepsilon\leqq 1$ such that 
the space $p_K\bigl(\pi(s)V^K\bigr)$ does not depend on 
$s\in S_I^-(\varepsilon)$ and is isomorphic to 
$(V_P)^{M_K}$ by the restriction of $j_P: V\to V_P$. 
The vectors in $p_K\bigl(\pi(s)V^K\bigr)$ are called canonical liftings 
over $(V_P)^{M_K}$ with respect to $K$. First, using 
(1) of Lemma \ref{1.6}, we can even construct 
the space $p_K\bigl(\pi(s)V^K\bigr)$ of 
canonical liftings by taking 
$s$ from $S_{0,I}^-(\varepsilon)$ (replacing $\varepsilon$ suitably). 
This step is entirely the same as the derivation of \cite[4.3.2]{C} from 
\cite[4.3.1]{C}. Next, choose $K$ small enough so that $V_1\subset V^K$. 
Then, $\pi(s)v\in V^{M_{0,K}U_{0,K}^-}$ for each $v\in V_1$ and 
$s\in S_{0,I}^-(\varepsilon)$. 
By \cite[5.3~(1)]{KT}, we have 
$$
\langle \lambda, \pi(s)v\rangle 
=\langle \lambda, p_{U_{0,K}}\bigl(\pi(s)v\bigr)\rangle
=\langle \lambda, p_{K}\bigl(\pi(s)v\bigr)\rangle. 
$$
Since $p_{K}\bigl(\pi(s)v\bigr)$ is 
now a vector in the space of canonical liftings, this is further equal to 
\begin{equation} \label{3.3.1}
\langle r_P(\lambda), j_P\bigl(p_K(\pi(s)v)\bigr)\rangle
=\langle r_P(\lambda), j_P\bigl(p_{U_{0,K}}(\pi(s)v)\bigr)\rangle
\end{equation} 
by the definition \eqref{3.2.1} of $r_P(\lambda)$. We use the decomposition 
$U_0=U\cdot U'$ where 
$U'=M\cap U_0$. 
This implies that $U_{0,K}=U_K\cdot U'_K$, hence 
$p_{U_{0,K}}=p_{U_K}\circ p_{U'_K}$. Since 
$j_P\circ p_{U_K}=j_P$ and $U'_K\subset M$, the right hand side 
of \eqref{3.3.1} is equal to 
$$
\langle r_P(\lambda), j_P\bigl(p_{U'_K}(\pi(s)v)\bigr)\rangle
=\langle r_P(\lambda), p_{U'_K}\bigl(j_P(\pi(s)v)\bigr)\rangle
$$
$$
=\langle r_P(\lambda), p_{U'_K}\bigl(\delta_P^{1/2}(s)\pi_P(s)j_P(v))\bigr)\rangle. 
$$
Finally, as a vector in the representation $\pi_P$ of $M$, $j_P(v)$ is ${M_K}$-fixed 
and thus $\pi_P(s)j_P(v)\in (V_P)^{M_{0,K}(U_{0,K}^-\cap M)}$. Applying 
\cite[5.3~(1)]{KT} for the 
$M\cap H$-invariant linear form $r_P(\lambda)$, we have 
$$
\langle r_P(\lambda), p_{U'_K}\bigl(\pi_P(s)j_P(v))\bigr)\rangle
=\langle r_P(\lambda), \pi_P(s)j_P(v)\rangle. 
$$
This completes the proof. 
\end{proof} 

\begin{Remark} \label{3.4} 
When we take a $\sigma$-split parabolic subgroup $P$ without 
specifying the initial data $({\bf S}_0, {\bf A}_{\emptyset}, 
{\Delta})$, we may think of $P$ as a 
standard one $P_{{I}}$ 
for a suitable choice of  $({\bf S}_0, {\bf A}_{\emptyset}, 
{\Delta})$ and $I\subset\Delta$. However, sometimes 
we have to fix the data 
$({\bf S}_0, {\bf A}_{\emptyset}, 
{\Delta})$ and deal with an arbitrary (possibly non-standard) 
$\sigma$-split parabolic subgroup written in the form 
$P=\gamma^{-1} P_{{I}}\gamma$ with 
$\gamma\in ({\bf M}_0{\bf H})(F)$. 
The relation as in \ref{3.3} for such a $P$ can be 
derived in the same way, and is presented as follows: There exists 
a positive real number $\varepsilon
=\varepsilon_{\gamma, I}\leqq 1$ such 
that 
\begin{equation}\label{3.4.1}
\langle\lambda,\pi(\gamma^{-1} s\gamma)v\rangle
=
\delta_P^{1/2}(\gamma^{-1} s\gamma)
\langle r_P(\lambda),\pi_P(\gamma^{-1}s\gamma)j_P(v)\rangle
\end{equation} 
for all $s\in S_{0, I}^-(\varepsilon)$, $v\in V_1$, and $\lambda\in(V^*)^H$. 
\end{Remark} 

\section{The main theorem} 

In this section we prepare notation on exponents, 
give a preliminary result (Proposition \ref{4.3}) 
on the relation between $H$-square integrability and 
exponents, and establish a 
criterion for $H$-square integrability (Theorem \ref{4.7}). 
We also give a non-trivial relation between 
$H$-square integrability and the usual square integrability 
(Proposition \ref{4.10}). 
 
\subsection{Exponents} \label{4.1} 

\indent 

Let $Z_1$ be a closed subgroup of the center of $G$ and 
$\mathcal X(Z_1)$ be the set of all quasi-characters of $Z_1$. 
For a smooth representation $(\pi, V)$ 
of $G$ and a quasi-character $\omega\in\mathcal X(Z_1)$, we put 
$$
V_{\omega,\infty}=
\left\{ v\in V\left|\,\begin{matrix} \text{There exists a }d\in\Bbb N \text{ such that}\\ 
\bigl(\pi(z)-\omega(z)\bigr)^d v=0\text{ for all }z\in Z_1
\end{matrix}\right. \right\}. 
$$
This is a $G$-submodule of $V$. 
We put 
$$
{\mathcal Exp}_{Z_1}(V)={\mathcal Exp}_{Z_1}(\pi)=
\left\{\omega\in{\mathcal X}(Z_1)\bigm| 
V_{\omega,\infty}\neq 0\right\}. 
$$
If $(\pi,V)$ is finitely generated and admissible, then the set 
${\mathcal Exp}_{Z_1}(\pi)$ is finite and $V$ has 
a direct sum decomposition 
$$
V=\bigoplus_{\omega\in{\mathcal Exp}_{Z_1}(\pi)}
V_{\omega,\infty}
$$
(see \cite[2.1.9]{C}). Let $Z_1$ and $Z_2$ be closed subgroups of the 
center of $G$ such that $Z_1\supset Z_2$. As is easily seen, 
the mapping 
$$ 
{\mathcal Exp}_{Z_1}(\pi)\to {\mathcal Exp}_{Z_2}(\pi)
$$ 
defined by restriction is surjective. 

\subsection{Exponents along parabolic subgroups} \label{4.2} 

\indent 

Let $(\pi,V)$ be a finitely generated admissible representation of $G$. 
For each parabolic subgroup $P$ of $G$ with the 
$F$-split component $A$, 
we consider the finite set 
${\mathcal Exp}_{A}(\pi_P)$. 
The elements of ${\mathcal Exp}_{A}(\pi_P)$ are called 
{\it exponents of $\pi$ along $P$}. 
If $P$ is a $\sigma$-split parabolic subgroup with the 
$(\sigma,F)$-split component $S$, we also consider the 
finite set ${\mathcal Exp}_{S}(\pi_{P})$. The mapping 
\begin{equation}\label{4.2.1} 
{\mathcal Exp}_{A}(\pi_{P})\to 
{\mathcal Exp}_{S}(\pi_{P})
\end{equation} 
defined by restriction is surjective.  

Let $P_1=M_1U_1$ and 
$P_2=M_2U_2$ be $\sigma$-split parabolic subgroups of $G$ 
with the $(\sigma,F)$-split components $S_1$ and $S_2$ 
respectively. If $P_1\subset P_2$, then 
$M_1\subset M_2$ and $S_1\supset S_2$. 
The intersection $M_2\cap P_1$ is a 
$\sigma$-split parabolic subgroup of $M_2$ 
having $M_1$ as a $\sigma$-stable 
Levi subgroup. As is well-known, 
$(\pi_{P_2})_{M_2\cap P_1}$ is 
naturally isomorphic to $\pi_{P_1}$ as an $M_1$-module. 
Through this isomorphism, it is easy to see that 
\begin{equation} \label{4.2.2}  
\chi \in {\mathcal Exp}_{S_1}(\pi_{P_1})\,\Longrightarrow \,
\chi |_{S_2} \in {\mathcal Exp}_{S_2}(\pi_{P_2}). 
\end{equation}

In the course of the proofs in this section, 
we need to fix the initial data 
$({\bf S}_0, {\bf A}_{\emptyset}, {\Delta})$ as in \ref{1.4}. 
In dealing with an arbitrary $\sigma$-split 
parabolic subgroup $P$, 
we have to 
write it as $P=\gamma^{-1}P_I\gamma$ by a 
$\sigma$-split subset $I\subset\Delta$ and 
an element $\gamma\in({\bf M}_0{\bf H})(F)$. 
The $(\sigma,F)$-split component of 
$P$ is then written as $S=\gamma^{-1}S_I\gamma$. 
We also write $S^-=\gamma^{-1}S_I^-(1)\gamma$ 
in such a case. 

Now, for a given 
finitely generated admissible 
representation $(\pi,V)$ of $G$, let us 
consider the following condition 
on a $\sigma$-split parabolic subgroup $P$: 
\begin{equation}\label{sharp1}
|\chi(s)|<1\text{ for all }\chi\in{\mathcal Exp}_{S}(\pi_P)
\text{ and all } 
s\in S^-\setminus Z_0S^1. \tag{$\sharp_P$} 
\end{equation} 

\begin{Proposition} \label{4.3} 
Let $\omega_0$ be a unitary character of $Z_0$ and 
$(\pi,V)$ a finitely generated $H$-distinguished 
admissible $\omega_0$-representation of $G$. 
If the condition \eqref{sharp1} is satisfied for every 
$\sigma$-split parabolic subgroup 
$P$ of $G$, then 
$(\pi, V)$ is $H$-square integrable with respect to any 
$\lambda\in\bigl(V^*\bigr)^H$. 
\end{Proposition} 
\begin{proof}  
Let $\Gamma$, $K_0$, $\Xi$ be as in 
\ref{2.2} for the data $({\bf S}_0, {\bf A}_{\emptyset}, {\Delta})$. 
Take a non-zero vector $v_0\in V$ and 
let $V_0$ be the finite dimensional 
subspace of $V$ generated by $\pi(k^{-1}\xi^{-1})v_0$ ($k\in K_0$, 
$\xi\in\Xi$). Further, let $V_1$ be the finite dimensional 
subspace of $V$ generated by $\pi(\gamma^{-1})v$ 
($\gamma\in\Gamma$, $v\in V_0$). 
We take a positive real number $\varepsilon\leqq 1$ satisfying 
\begin{itemize} 
\item $\varepsilon\leqq \varepsilon_{\gamma, I}$ for all $\gamma\in\Gamma$ and 
all $\sigma$-split $I\subset \Delta$ 
where $\varepsilon_{\gamma,I}$ is such that 
\eqref{3.4.1} is valid for all $v\in V_1$, and 
\item $\varepsilon\leqq c_{{I}}$ for all $\sigma$-split 
${I}\subset {\Delta}$ 
where $c_{{I}}$ is the constant given in the proof of Lemma \ref{1.6}. 
\end{itemize} 
We use a disjoint decomposition 
$$
S_0^+=\bigcup_{I\subset{\Delta}:\sigma\text{-split}} {}_IS_0^+(\varepsilon)
$$
obtained from Lemma \ref{1.6} (3) for the number $\varepsilon$ as above. 
Let us put 
$$
G_{I,\gamma}=\bigcup_
{\dot{s}\in {}_IS_{0}^+(\varepsilon)/Z_0S_0^1} 
K_0 \dot{s}\gamma H
$$
for each $\sigma$-split subset $I\subset{\Delta}$ 
and $\gamma\in\Gamma$. Then, 
by \eqref{2.2.1} we have 
$$
G/Z_0H=
\bigcup_{\xi,\gamma, {I}} \xi G_{I,\gamma}/H. 
$$
Now we start to evaluate the $L^2$-norm of the $H$-matrix coefficient 
$\varphi_{\lambda, v_0}$. It is clear that 
$$
\int_{G/Z_0H}\bigl|\varphi_{\lambda,v_0}(g)\bigr|^2dg
\leqq 
\sum_{\xi,\gamma, {I}}
\left(
\int_{\xi G_{{I},\gamma}/H}\bigl|\varphi_{\lambda,v_0}(g)\bigr|^2dg\right). 
$$
It is enough to study the convergence of 
$$
\int_{\xi G_{{I},\gamma}/H}\bigl|\varphi_{\lambda,v_0}(g)\bigr|^2dg
=\int_{G_{{I},\gamma}/H}\bigl|\varphi_{\lambda,\pi(\xi^{-1})v_0}(g)\bigr|^2dg
$$
for each $\xi$, $\gamma$ and ${I}$. So the proof of the proposition 
is completed once the following claim is shown. 
\begin{Claim} 
If \eqref{sharp1} is satisfied for $P=\gamma^{-1} P_{{I}}\gamma$, 
then 
$$
\int_{G_{I,\gamma}/H}\bigl|\varphi_{\lambda,v}(g)\bigr|^2dg<\infty
$$ 
for all $v\in V_0$. 
\end{Claim} 
Let us prove this. It is obvious that 
\begin{equation} \label{4.3.1} 
\int_{G_{I,\gamma}/H}\bigl|\varphi_{\lambda,v}(g)\bigr|^2dg
\leqq 
\sum_{\dot{s}\in
{}_I S_0^+(\varepsilon)/Z_0S_0^1}
\int_{K_0\dot{s}\gamma H/H}\bigl|\varphi_{\lambda,v}(g)\bigr|^2dg. 
\end{equation} 
There is an element 
$k_0\in K_0$ such that 
$$
\bigl|\varphi_{\lambda,v}(k\dot{s}\gamma h)\bigr|^2\quad(k\in K_0) 
$$
attains the maximum at $k=k_0$. We have 
\begin{align} \label{4.3.2}
\int_{K_0\dot{s}\gamma H/H}\bigl|\varphi_{\lambda,v}(g)\bigr|^2dg
&\leqq {\rm{vol}}(K_0\dot{s}\gamma H/H)\cdot 
\bigl|\varphi_{\lambda,v}(k_0\dot{s}\gamma h)\bigr|^2 \\
&\leqq C\cdot\delta_{P_0}(\dot{s})\cdot 
\bigl|\varphi_{\lambda,\pi(k_0^{-1})v}(\dot{s}\gamma)\bigr|^2 \notag 
\end{align}
for some positive real constant $C$ which does not 
depend on $\dot{s}$ by \ref{2.6}. We may replace 
$\pi(k_0^{-1})v\in V_0$ by $v$. 
Applying Proposition \ref{3.3} (or \ref{3.4}) 
along $P=\gamma^{-1}P_{{I}}\gamma$, 
we have 
\begin{align}
\varphi_{\lambda,v}(\dot{s}\gamma)
&=\langle \lambda, \pi(\gamma^{-1}\dot{s}^{-1})v\rangle=
\langle \lambda, \pi(\gamma^{-1}\dot{s}^{-1}\gamma)\pi(\gamma^{-1})v\rangle \notag\\
&=\delta_P(\gamma^{-1}\dot{s}^{-1}\gamma)^{1/2}
\langle r_P(\lambda), \pi_P(\gamma^{-1}\dot{s}^{-1}\gamma)j_P
(\pi(\gamma^{-1})v)\rangle \notag
\end{align}
since $\dot{s}^{-1}\in{}_I S_0^-(\varepsilon)\subset S_{0,I}^-(\varepsilon)$ and 
$\pi(\gamma^{-1})v\in V_1$. 
Next, we use Lemma \ref{1.6} (2) to write 
$$
\dot{s}^{-1}=\dot{s}_I\cdot t_i,\quad \dot{s}_I\in S_I^-(\varepsilon')/Z_0S_I^1, 
$$
where $\varepsilon'\leqq 1$ by our choice of $\varepsilon$. Putting 
$\pi_P(\gamma^{-1}t_i\gamma)j_P
(\pi(\gamma^{-1})v)=\overline{v}$ for simplicity, 
we have 
$$
\varphi_{\lambda,v}(\dot{s}\gamma)
=\delta_{P_{{I}}}(\dot{s}_It_i)^{-1/2}
\langle r_P(\lambda), \pi_P(\gamma^{-1}\dot{s}_I \gamma)\overline{v}
\rangle. 
$$
The function 
$\langle r_P(\lambda), \pi_P(\cdot)\overline{v}
\rangle$ on $S=\gamma^{-1}S_I\gamma$ is $S$-finite. Thus it 
can be written as 
\begin{equation} \label{4.3.3}
\langle r_P(\lambda), \pi_P(s')\overline{v}
\rangle=\sum_{\chi\in{\mathcal Exp}_S(\pi_P)}
\chi(s'){\mathcal P}_{\chi}\bigl(\nu_{S}(s')\bigr)
\end{equation} 
for all $s'\in S$ using suitable polynomials ${\mathcal P}_{\chi}$ on 
$\frak a_{\bf S}$ (see \cite[I.2]{W}). Let us write 
$s'_I=\gamma^{-1}\dot{s}_I\gamma$. 
Then, returning to \eqref{4.3.2}, 
we have a bound for $\int_{K_0\dot{s}\gamma H/H}
\bigl|\varphi_{\lambda,v}(g)\bigr|^2dg$ by 
\begin{align}
C&\cdot 
\delta_{P_0}(\dot{s}_I t_i)\delta_{P_{{I}}}(\dot{s}_It_i)^{-1}
\Bigl|\sum_{\chi}\chi(s'_I)
{\mathcal P}_{\chi}\bigl(\nu_{S}(s'_I)\bigr)\Bigr|^2 \notag\\
&=C\cdot \delta_{P_0}(t_i)\delta_{P_{{I}}}(t_i)^{-1}
\Bigl|\sum_{\chi}\chi(s'_I)
{\mathcal P}_{\chi}\bigl(\nu_{S}(s'_I)\bigr)\Bigr|^2\notag\\ 
&\leqq 
C\cdot \delta_{P_0}\delta_{P_{{I}}}^{-1}(t_i)
\sum_{\chi,\chi'}\Bigl|\chi\chi'(s'_I){\mathcal P}_{\chi}\bigl({\nu}_{S}(s'_I)
\bigr){\mathcal P}_{\chi'}\bigl({\nu}_{S}(s'_I)
\bigr)\Bigr|.\notag
\end{align} 
Here we used the fact that $\delta_{P_0}\equiv\delta_{P_{{I}}}$ on 
$S_I$. Returning further to \eqref{4.3.1}, the 
integral 
$\int_{G_{I,\gamma}/H}\bigl|
\varphi_{\lambda,v}(g)\bigr|^2dg$ is bounded by 
$$
\sum_{\dot{s}_I\in S_I^-(\varepsilon')/Z_0S_I^1} 
\sum_i 
\delta_{P_0}\delta_{P_{{I}}}^{-1}(t_i)
\sum_{\chi,\chi'}\Bigl|\chi\chi'(s'_I)
{\mathcal P}_{\chi}\bigl({\nu}_{S}(s'_I)
\bigr){\mathcal P}_{\chi'}\bigl({\nu}_{S}(s'_I)\bigr)\Bigr|. 
$$
Now, we may regard $S_I^-(\varepsilon')/Z_0S_I^1$ as 
a set of lattice points in a positive cone. Then the infinite sum 
$$
\sum_{\dot{s}_I\in S_I^-(\varepsilon')/Z_0S_I^1}
\Bigl|\chi\chi'(s'_I){\mathcal P}_{\chi}\bigl({\nu}_{S}(s'_I)
\bigr){\mathcal P}_{\chi'}\bigl({\nu}_{S}(s'_I)\bigr)\Bigr|
$$
is essentially a power series with polynomial coefficients. 
This series converges if $\bigl|\chi(s'_I)\bigr|<1$ for 
all $\chi$ and $s'_I=\gamma^{-1}\dot{s}_I\gamma$ with 
$\dot{s}_I\in S_I^-(\varepsilon')/Z_0S_I^1$, 
except for those $\dot{s}_I\in Z_0S_I^1$ 
which represent the origin in the lattice. Thus the condition 
\eqref{sharp1} for $P=\gamma^{-1}P_I\gamma$ is sufficient 
for the convergence of 
$\int_{G_{I,\gamma}/H}\bigl|\varphi_{\lambda,v}(g)\bigr|^2dg$. 
\end{proof} 

\subsection{Exponents with respect to $\lambda$}%{4.4} 

\indent 

Let $(\pi,V)$ be a smooth representation of 
$G$ and take an $H$-invariant linear form $\lambda$ 
on $V$. We define 
$$
{\mathcal Exp}_{Z_0}(\pi, \lambda)
:=\left\{ \omega\in {\mathcal Exp}_{Z_0}(\pi)\Bigm| 
\lambda|_{V_{\omega,\infty}}\neq 0\right\}. 
$$ 
Assume that $(\pi,V)$ is finitely generated and admissible. 
For each $\sigma$-split parabolic subgroup 
$P$ of $G$ with the 
$(\sigma,F)$-split component $S$, 
we consider the subset 
${\mathcal Exp}_S\bigl(\pi_P, r_P(\lambda)\bigr)$ 
of ${\mathcal Exp}_S(\pi_P)$. As a relation similar to 
\eqref{4.2.2}, we have 
the following. 
\begin{Lemma} \label{4.5}  
Let $P_1$ and $P_2$ be $\sigma$-split parabolic subgroups of $G$ such that 
$P_1\subset P_2$, with the $(\sigma,F)$-split components $S_1$ and $S_2$ 
respectively. Then, 
$$
\chi\in {\mathcal Exp}_{S_1}\bigl(\pi_{P_1}, r_{P_1}(\lambda)\bigr)
\,\Longrightarrow \, 
\chi|_{S_2}\in
{\mathcal Exp}_{S_2}\bigl(\pi_{P_2}, r_{P_2}(\lambda)\bigr). 
$$
\end{Lemma} 
\begin{proof}  
Suppose that $\chi\in {\mathcal Exp}_{S_1}\bigl(\pi_{P_1}, r_{P_1}(\lambda)\bigr)$. 
For each $\bar{v}_1\in (V_{P_1})_{\chi,\infty}$, we can take 
a $\bar{v}_2\in (V_{P_2})_{\chi|_{S_2},\infty}$ 
such that $j_{M_2\cap P_1}(\bar{v}_2)=\bar{v}_1$ 
(in the notation of \ref{4.2}) regarding \eqref{4.2.2}. 
Applying Proposition \ref{3.3} to 
the $M_2\cap H$-matrix coefficients, we have 
$$
\langle r_{P_2}(\lambda), \pi_{P_2}(s)\bar{v}_2\rangle 
=\delta_{M_2\cap P_1}(s)^{1/2} 
\langle r_{M_2\cap P_1}\bigl(r_{P_2}(\lambda)\bigr), \pi_{P_1}(s)\bar{v}_1\rangle 
$$
at least for some $s\in S_1$. The right hand side is 
written as 
$$
\delta_{M_2\cap P_1}(s)^{1/2} 
\langle r_{P_1}(\lambda), \pi_{P_1}(s)\bar{v}_1\rangle 
$$
by the transitivity result $r_{M_2\cap P_1}\circ r_{P_2}=r_{P_1}$ given in 
\cite[Proposition 5.9]{KT} (and also \cite[Th\'{e}or\`{e}me 3]{L}). 
Hence $r_{P_2}(\lambda)$ cannot be 
identically zero on ${(V_{P_2})_{\chi|_{S_2},\infty}}$ 
provided that 
$r_{P_1}(\lambda)|_{(V_{P_1})_{\chi,\infty}}\neq 0$. 
\end{proof} 

Fix a non-zero $\lambda\in (V^*)^H$. 
Let us take up the following condition on 
a $\sigma$-split parabolic subgroup $P$: 
\begin{equation} \label{sharp2} 
|\chi(s)|<1\text{ for all }\chi\in{\mathcal Exp}_S\bigl(\pi_P, r_P(\lambda)\bigr) 
\text{ and all } 
s\in S^-\setminus Z_0S^1. 
\tag{$\sharp_{P,\lambda}$} 
\end{equation}

\begin{Lemma} \label{4.6} 
If the condition \eqref{sharp2} holds for 
every maximal $\sigma$-split parabolic subgroup $P$ of $G$, then 
it holds for every $\sigma$-split parabolic subgroup $P$ of $G$. 
\end{Lemma} 
\begin{proof}  
It is enough to derive \eqref{sharp2} for a standard 
$P=P_I$ assuming \eqref{sharp2} for all maximal standard 
$P$. For a given 
$\sigma$-split subset $I$ of $\Delta$, let us write 
$$
\overline{\Delta}\setminus\overline{I}
=\{\overline{\alpha}_1,\,\overline{\alpha}_2,\,\dots,\,\overline{\alpha}_r\}. 
$$
For each $i$, let ${\bf P}_i$ be the maximal standard 
$\sigma$-split parabolic subgroup corresponding to 
$[\overline{\Delta}\setminus\{\overline{\alpha}_i\}]$ (see the remark 
at the end of \ref{1.4}). 
These are maximal $\sigma$-split parabolic 
subgroups containing ${\bf P}_I$. 
The $(\sigma,F)$-split component ${\bf S}_i$ of 
${\bf P}_i$ is 
given by 
\begin{equation} \label{4.6.1} 
{\bf S}_i={\bf S}_{\overline{\Delta}\setminus 
\{\overline{\alpha}_i\}}=
\left(\bigcap_{\overline{\alpha}\neq\overline{\alpha}_i} \ker (\overline{\alpha})
\right)^0. 
\end{equation} 
The $(\sigma,F)$-split component ${\bf S}={\bf S}_I$ of ${\bf P}={\bf P}_I$ can be 
decomposed as ${\bf S}={\bf S}_1{\bf S}_2\dots{\bf S}_r$, so that 
$S_1S_2\dots S_r$ is of finite index 
in $S$. Thus for any given $s\in S^-$, there exists an integer $m$ such that 
$s^m\in S_1S_2\dots S_r$. We may write 
$s^m=s_1s_2\dots s_r$ ($s_i\in S_i$). From \eqref{4.6.1} it is 
easy to see that $s_i\in S_i^-$ for each $i$. Now, for any 
$\chi\in{\mathcal Exp}_{S}\bigl(\pi_{P}, r_{P}(\lambda)\bigr)$, 
we have 
$\chi|_{S_i}\in{\mathcal Exp}_{S_i}\bigl(\pi_{P_i}, r_{P_i}(\lambda)\bigr)$ 
by Lemma \ref{4.5}. Therefore, 
assumptions \eqref{sharp2} for all $P=P_i$ imply that 
$$
\bigl|\chi(s)\bigr|^m=\bigl|\chi(s^m)\bigr|
=\bigl|\chi(s_1)\bigr|\bigl|\chi(s_2)\bigr|\dots
\bigl|\chi(s_r)\bigr|< 1. 
$$
This completes the proof. 
\end{proof} 

Now we shall state the main theorem of this 
paper. 

\begin{Theorem} \label{4.7} 
Let $\omega_0$ be a 
unitary character of $Z_0$ and 
$(\pi,V)$ a finitely generated $H$-distinguished admissible 
$\omega_0$-representation of $G$. 
Then, for a non-zero $H$-invariant linear form 
$\lambda$ on $V$, the representation $(\pi,V)$ is 
$H$-square integrable with respect to $\lambda$ if and only if 
the condition \eqref{sharp2} is 
satisfied for every $\sigma$-split parabolic subgroup 
$P$ of $G$. 
\end{Theorem} 
\begin{proof}  
Note that only the exponents in 
${\mathcal Exp}_S\bigl(\pi_P, r_P(\lambda)\bigr)$ 
contribute to the evaluation of 
$\int_{G/Z_0H} |\varphi_{\lambda,v}(g)|^2dg$ 
in the proof of Proposition \ref{4.3}, 
specifically at the expression \eqref{4.3.3}. 
So the {\it if part} is already proved in \ref{4.3}. 
To prove the {\it only if part}, suppose that \eqref{sharp2} 
fails for some $\sigma$-split parabolic subgroup $P$. 
By Lemma \ref{4.6}, we may suppose that $P$ 
is a maximal one, say, $P=\gamma^{-1}P_I\gamma$ with 
$I=[\overline{\Delta}\setminus\{\overline{\alpha}\}]$ for some 
$\overline{\alpha}\in\overline{\Delta}$. 
Then there exists an exponent $\chi\in
{\mathcal Exp}_{S}\bigl(\pi_{P}, r_{P}(\lambda)\bigr)$ and an 
element $s=\gamma^{-1}s_I\gamma\in S^-\setminus Z_0S^1$, 
$s_I\in S_I^-\setminus Z_0S_I^1$, 
such that $\bigl|\chi(s)\bigr|\geqq 1$. 
Take a vector $\bar{v}\in(V_P)_{\chi,\infty}$ with 
$\langle r_P(\lambda), \bar{v}\rangle\neq 0$ and let 
$v\in V$ be such that $j_P(\pi(\gamma^{-1})v)=\bar{v}$. 

From Proposition \ref{2.3}, the union 
$\bigcup_{n\geqq 0} K_0s_I^{-n}\gamma H$ is disjoint. So it is enough 
to see that 
$$
\sum_{n\geqq 0} \int_{K_0 s_I^{-n}\gamma H/H} 
\bigl| \varphi_{\lambda, v}(g)\bigr|^2dg
$$
is divergent. Take an open compact subgroup 
$K\subset K_0$ which fixes $v$. 
For each $n$, the function $\varphi_{\lambda,v}$ 
is constant on $Ks_I^{-n}\gamma H$, hence 
\begin{align} \label{4.7.1} 
\int_{K_0 s_I^{-n}\gamma H/H} 
\bigl| \varphi_{\lambda, v}(g)\bigr|^2dg
&\geqq \int_{K s_I^{-n}\gamma H/H} 
\bigl| \varphi_{\lambda, v}(g)\bigr|^2dg \notag \\ 
&={\rm vol}(Ks_I^{-n}\gamma H/H)\bigl| \langle 
\lambda, \pi(\gamma^{-1}s_I^n)v\rangle \bigr|^2 \notag \\
&=C\cdot\delta_{P_0}(s_I^{-n})\bigl| \langle 
\lambda, \pi(\gamma^{-1}s_I^n\gamma)\pi(\gamma^{-1})v\rangle \bigr|^2 
\end{align} 
for some constant $C$ by Proposition \ref{2.6}. 
Now, let us take $\varepsilon\leqq 1$ such that 
the relation \eqref{3.4.1} holds for all $s\in S_I^-(\varepsilon)$. Since $I=[\overline{\Delta}\setminus\{\overline{\alpha}\}]$ here, 
the set $S_I^-(\varepsilon)$ is described as 
$$
S_I^-(\varepsilon)=\left\{ s\in S_I\,\left| \,
\bigl| s^{\overline{\alpha}}\bigr|_F\leqq \varepsilon\right.\right\}. 
$$
The element 
$s_I\in S^-_I\setminus Z_0S_I^1$ satisfies $
\bigl| (s_I)^{\overline{\alpha}}\bigr|_F<1$. Thus we can take an 
integer $N$ such that $s_I^n\in S_I^-(\varepsilon)$ for all 
$n\geqq N$. As a result, \eqref{4.7.1} continues as 
\begin{align} 
C\cdot\delta_{P_0}(s_I^{-n})&\delta_{P_I}(s_I^n)
\bigl| \langle r_P(\lambda), 
\pi_P(\gamma^{-1}s_I^n\gamma)j_P(\pi(\gamma^{-1})v)\rangle \bigr|^2 
\notag \\
&=C\cdot \bigl| \langle r_P(\lambda), 
\pi_P(s^n)\bar{v}\rangle \bigr|^2 \notag
\end{align}
whenever $n\geqq N$. Since $\bar{v}\in(V_P)_{\chi,\infty}$, the function 
$\langle r_P(\lambda), \pi_P(\cdot)\bar{v}\rangle$ on $S$ 
is written just as $\chi(\cdot){\mathcal P}_{\chi}\bigl(\nu_S(\cdot)\bigr)$ by a 
single polynomial ${\mathcal P}_{\chi}$ on $\frak a_{\bf S}$. Finally, the series
$$
\sum_{n\geqq N} \bigl|\chi(s^n){\mathcal P}_{\chi}(\nu_S(s^n))
\bigr|^2=\sum_{n\geqq N} \bigl|\chi(s)\bigr|^{2n}
\bigl|{\mathcal P}_{\chi}(\nu_S(s^n))
\bigr|^2
$$
is obviously divergent, hence the proof is completed. 
\end{proof} 

\begin{Remark} %{4.8} 
We have actually shown that the following 
three conditions are equivalent: 
\begin{enumerate} {\it 
\item $(\pi, V)$ is $H$-square integrable with respect to 
$\lambda$. 
\item \eqref{sharp2} is 
satisfied for every $\sigma$-split parabolic subgroup 
$P$ of $G$. 
\item \eqref{sharp2} is 
satisfied for every maximal $\sigma$-split parabolic subgroup 
$P$ of $G$. }
\end{enumerate} 
See \cite[III.1.1]{W} for a similar statement 
on the usual square integrability. 
\end{Remark} 

\subsection{$H$-square integrability and the usual square integrability}%{4.9} 

\indent 

Let $(\widetilde{\pi},\widetilde{V})$ be the contragredient of 
$(\pi,V)$. For $v\in V$ and 
$\widetilde{v}\in\widetilde{V}$, the usual 
matrix coefficient $c_{\widetilde{v},v}$ is 
defined by 
$$
c_{\widetilde{v},v}(g)=\langle \widetilde{v},\pi(g^{-1})v\rangle.  
$$
Let $\omega$ be a unitary character of the $F$-split component $Z$ 
of $G$. A smooth $\omega$-representation $(\pi,V)$ of $G$ is 
said to be square integrable 
if the functions $\bigl|c_{\widetilde{v},v}(\cdot)\bigr|$ are 
square integrable on $G/Z$ for all $v$ and $\widetilde{v}$. 
For the $F$-split component 
$A=A_I$ of a parabolic subgroup 
$P=P_I$ of $G$, the dominant 
part $A^-$ of $A$ is defined by 
$$
A^-=\left\{ a\in A\bigm| \left|a^{\alpha}\right| \leqq 1 \,
(\alpha\in\Delta\setminus I)\right\}.  
$$
There is a well-known criterion for square integrability 
due to Casselman (\cite[4.4.6]{C}). In terms of normalized 
Jacquet modules, it is stated as follows. 

\begin{Cass} \label{Cass}  
A finitely generated admissible $\omega$-representation 
$(\pi,V)$ of $G$ is square integrable if and 
only if for every parabolic subgroup $P$, 
$\bigl|\chi(a)\bigr|<1$ holds 
for all $\chi\in{\mathcal Exp}_{A}(\pi_P)$ and 
$a\in A^-\setminus ZA^1$. 
\end{Cass} 

Although the usual square integrability and $H$-square integrability are 
different notions, we have the following relation.  
 
\begin{Proposition} \label{4.10} 
Let $\omega$ be a unitary character of $Z$ and 
$(\pi,V)$ a finitely generated admissible $\omega$-representation 
of $G$. If $(\pi,V)$ is square integrable and 
is $H$-distinguished, 
then it is $H$-square integrable 
for all $\lambda\in\bigl(V^*\bigr)^H$. 
\end{Proposition} 

\begin{proof}  
This is immediate from Casselman's criterion and 
Proposition \ref{4.3}, 
in view of the surjectivity of \eqref{4.2.1}. 
\end{proof} 

An example of square integrable $H$-distinguished representations 
can be found in \cite[\S7]{H} for the 
symmetric space 
${\rm GL}_2(E)/{\rm GL}_2(F)$ where $E/F$ is a quadratic extension. 

\section{Examples of $H$-square integrable representations} 

We shall give two simple examples of $H$-square integrable 
representations which are 
not square integrable. In both cases, 
there are constructions (\cite{D} and \cite{HR}) 
of an $H$-distinguished representation $\pi(\rho)$ 
attached to a representation $\rho$ of ${\rm GL}_2(F)$. 
We have observed in \cite[\S 8]{KT} that 
$\pi(\rho)$ is $H$-relatively cuspidal if 
$\rho$ is cuspidal. In what follows, we shall 
observe that 
$\pi(\rho)$ is $H$-square integrable if 
$\rho$ is square integrable. 

\subsection{The symmetric space 
${\rm GL}_3(F)/\bigl({\rm GL}_2(F)\times {\rm GL}_1(F)\bigr)$} \label{5.1} 

\indent 

Let $G$ be the group ${\rm GL}_3(F)$ and $\sigma$ the inner involution 
$\sigma={\rm{Int}}(\epsilon)$ 
defined by the anti-diagonal permutation matrix $\epsilon
=\left(\begin{smallmatrix} & &1\\ &1& \\ 1& & \end{smallmatrix}\right)$. 
Then the $\sigma$-fixed point subgroup $H$ 
is isomorphic to ${\rm GL}_2(F)\times {\rm GL}_1(F)$. 
For this symmetric space, 
all the irreducible $H$-distinguished representations were 
determined by D. Prasad \cite{P}. We shall use notation in \cite[8.2]{KT} for the case $n=3$. 

Let $Q=LU_Q$ be the standard parabolic subgroup of $G$ of type $(1,2)$. 
Thus $L\simeq{\rm GL}_{1}(F)\times{\rm GL}_{2}(F)$. 
Let $\rho$ be an infinite dimensional irreducible 
admissible representation of ${\rm GL}_2(F)$ with 
trivial central character and form the 
normalized induction 
$$
\pi({\rho})={\rm Ind}_{Q}^G (1_{{\rm GL}_{1}(F)}\otimes \rho) 
$$
where $1=1_{{\rm GL}_{1}(F)}$ denotes the trivial character of ${\rm GL}_{1}(F)
=F^{\times}$. 
Then $\pi({\rho})$ is irreducible and $H$-distinguished (\cite[Theorem 2 (2)]{P}). 

We take a maximal $(\sigma,F)$-split torus $S_0$ as the 
one consisting of diagonal matrices of the form 
${\rm diag}(s, 1, s^{-1})$ with $s\in F^{\times}$, and 
a maximal $F$-split torus $A_{\emptyset}$ of all diagonal matrices. 
The $(\sigma,F)$-split component $Z_0$ of $G$ is trivial. 
As a minimal 
$\sigma$-split parabolic subgroup $P_0$, we may take the 
Borel subgroup consisting of upper triangular matrices. 
In this case, $P_0$ is the only proper $\sigma$-split parabolic 
of $G$ up to $H$-conjugacy. Let us determine the 
exponents of $\pi(\rho)$ along $P_0$. 

By the Geometric Lemma of \cite[2.12]{BZ}, we have 
$$
\bigl(\pi(\rho)\bigr)_{P_0}^{\rm s.s}=
\left({\rm Ind}_{Q}^G (1\otimes \rho)\right)_{P_0}^{\rm s.s} 
\simeq 
\bigoplus_{w\in[W_L\backslash W]} 
\mathcal F_w(1\otimes \rho)^{\rm s.s}
$$
where $W$ and $W_L$ denote the Weyl group of $A_{\emptyset}$ in $G$ and 
$L$ respectively, $[W_L\backslash W]$ as in \cite[1.1.3]{C}, $\bigl(\cdots\bigr)^{\rm s.s}$ 
denotes the semisimplified form, and 
$$
\mathcal F_w(1\otimes \rho)
={\rm Ind}_{M_0\cap wQw^{-1}}^{M_0} 
\bigl(w\cdot(1_{{\rm GL}_{1}(F)}\otimes \rho)_{L\cap w^{-1}P_0w}\bigr). 
$$
The set $[W_L\backslash W]$ consists of three elements. 
In forms of permutation matrices, those are 
$$
e=\left(\begin{smallmatrix} 1& & \\ &1& \\
 & &1\end{smallmatrix}\right), \quad 
\left(\begin{smallmatrix} &1& \\ 1& & \\
 & &1\end{smallmatrix}\right),\quad 
\text{and }
\left(\begin{smallmatrix} &1& \\ & &1 \\
1& & \end{smallmatrix}\right). 
$$
It is easy to compute the $M_0=A_{\emptyset}$-module 
$\mathcal F_w(1\otimes \rho)$ for each $w$: 
$$
\mathcal F_w(1\otimes \rho)=
w\cdot\bigl(1\otimes \rho_{B_2}\bigr) 
$$
where $B_2$ denotes the Borel subgroup 
consisting of 
upper triangular matrices of ${\rm GL}_2(F)$. 

Now take $\rho$ to be the Steinberg representation 
${\rm St}_2$ of ${\rm GL}_2(F)$. We shall 
observe that 
{\it the irreducible $H$-distinguished 
representation 
$\pi({\rm St}_2)$ is 
$H$-square integrable.} 
As is well-known, 
$({\rm St}_2)_{B_2}$ is given by the 
one-dimensional character 
$\delta_{B_2}^{1/2}$. So the set 
${\mathcal Exp}_{A_{\emptyset}}\bigl(\pi({\rm St}_2)_{P_0}\bigr)$ 
of exponents of $\pi({\rm St}_2)$ along $P_0$ consists of 
three characters $\chi_1$, $\chi_2$, 
$\chi_3$ respectively given by 
$$
\chi_1\left(\begin{smallmatrix} a_1& & \\ &a_2& \\
 & &a_3\end{smallmatrix}\right)=\bigl|a_2\bigr|_F\bigl|a_3\bigr|_F^{-1}, 
\quad 
\chi_2\left(\begin{smallmatrix} a_1& & \\ &a_2& \\
 & &a_3\end{smallmatrix}\right)=\bigl|a_1\bigr|_F\bigl|a_3\bigr|_F^{-1},
 $$
 $$
 \chi_3\left(\begin{smallmatrix} a_1& & \\ &a_2& \\
 & &a_3\end{smallmatrix}\right)=\bigl|a_1\bigr|_F\bigl|a_2\bigr|_F^{-1}. 
$$
Finally, since $S_0^-$ (resp. $S_0^1$) 
consists of ${\rm diag}(s, 1, s^{-1})$ with $|s|_F\leqq 1$ (resp. $|s|_F=1$), 
we may conclude that the 
restriction to $S_0$ of each exponent $\chi_i$ satisfies 
$\bigl|\chi_i(s)\bigr|<1$ for all $s\in S_0^-\setminus S_0^1$. 
Hence the claim follows from Proposition \ref{4.3}.  

\subsection{The symmetric space 
${\rm GL}_4(F)/{\rm Sp}_2(F)$}  \label{5.2} 

\indent 

Let $G$ be the group ${\rm GL}_4(F)$ and $\sigma$ the involution 
on $G$ defined by
$$
\sigma(g)=\left(\begin{smallmatrix} 0&1& & \\ 
-1&0& & \\ & &0&1\\ & &-1&0
\end{smallmatrix}\right){}^tg^{-1} 
\left(\begin{smallmatrix} 0&1& & \\ 
-1&0& & \\ & &0&1\\ & &-1&0
\end{smallmatrix}\right)^{-1} \quad(g\in G).  
$$
Then the $\sigma$-fixed point subgroup $H$ 
is the symplectic group ${\rm Sp}_2(F)$. 
For this symmetric space, $H$-distinguished representations were 
studied by Heumos-Rallis \cite{HR}. See also \cite[8.3]{KT}. 

We take 
$$
S_0=\{{\rm diag}(s_1,s_1,s_2,s_2)\bigm| s_1,\,s_2\in F^{\times}\} 
$$
as a maximal $(\sigma,F)$-split torus of $G$. 
The $(\sigma,F)$-split component $Z_0$ of $G$ consists of 
all the scalar matrices of $G$ in this case. Let 
$P_0=M_0U_0$ be the standard parabolic subgroup of $G$ 
of type $(2,2)$. This is the only proper 
$\sigma$-split parabolic subgroup of $G$ 
up to $H$-conjugacy. Note that 
\begin{equation} \label{4.9.1} 
M_0\simeq {\rm GL}_2(F)\times {\rm GL}_2(F), \quad 
M_0\cap H\simeq {\rm SL}_2(F)\times {\rm SL}_2(F). 
\end{equation} 
Let $\rho$ be an irreducible admissible 
representation of $G_2:={\rm GL}_2(F)$ and let us form the 
normalized induction 
$$
{\rm I}(\rho)={\rm Ind}_{P_0}^G\bigl(\rho\cdot \bigl|\det(\cdot)\bigr|_F^{1/2}
\otimes \rho\cdot \bigl|\det(\cdot)\bigr|_F^{-1/2}\bigr). 
$$
Then ${\rm I}(\rho)$ is $H$-distinguished (\cite[11.1(a)]{HR}). 
If further $\rho$ is square integrable, 
then ${\rm I}(\rho)$ has the unique irreducible 
quotient $\pi(\rho)$ which also is $H$-distinguished (\cite[11.1(b)]{HR}). 

We shall investigate the 
$M_0$-module $\bigl({\rm I}(\rho)\bigr)_{P_0}$ and 
$M_0\cap H$-distinguished components therein. Note that 
irreducible $M_0\cap H$-distinguished $M_0$-modules have to be 
one dimensional according to \eqref{4.9.1}. 

Let $W$ (resp. $W_{M_0}$) be the Weyl group of 
$(G, A_{\emptyset})$ (resp. of $(M_0, A_{\emptyset})$). 
Following the definition of \cite[1.1.3]{C}, we can give 
the coset representatives $[W_{M_0}\backslash W/W_{M_0}]$ as 
$$
e,\quad \left(\begin{smallmatrix} &&1&0 \\ 
 & &0&1 \\ 1&0& & \\ 0&1& & 
\end{smallmatrix}\right),\quad \text{and }\,
\left(\begin{smallmatrix} 1& & & \\ 
 &0&1& \\ &1&0& \\ & & &1
\end{smallmatrix}\right). 
$$
We use the abbreviation 
$[\rho]=
\rho\cdot \bigl|\det(\cdot)\bigr|_F^{1/2}
\otimes \rho\cdot \bigl|\det(\cdot)\bigr|_F^{-1/2}$. 
The Geometric Lemma \cite[2.12]{BZ} asserts that 
$$
\bigl({\rm I}(\rho)\bigr)_{P_0}^{\rm s.s} 
\simeq \bigoplus_{w\in[W_{M_0}\backslash W/W_{M_0}]} 
\mathcal F_w([\rho])^{\rm s.s} 
$$
where 
$$
\mathcal F_w([\rho])=
{\rm Ind}_{M_0\cap wP_0w^{-1}}^{M_0} 
\bigl({}^w ([\rho]_{M_0\cap w^{-1}P_0 w})
\bigr).  
$$ 
It is easy to determine 
the three pieces $\mathcal F_w([\rho])$ for each $w$. 

\noindent 
(i) $w=e$: Nothing but 
$$
\mathcal F_w([\rho])=[\rho], 
$$
which is not $M_0\cap H$-distinguished unless 
$\rho$ is one dimensional. 

\noindent 
(ii) $w=\left(\begin{smallmatrix} &&1&0 \\ 
 & &0&1 \\ 1&0& & \\ 0&1& & 
\end{smallmatrix}\right)$: In this case, we have 
$$
M_0\cap wP_0w^{-1}=M_0\cap w^{-1}P_0w=M_0, 
$$
so that 
$$
\mathcal F_w([\rho])={}^w [\rho]= 
\rho\cdot \bigl|\det(\cdot)\bigr|^{-1/2}
\otimes \rho\cdot \bigl|\det(\cdot)\bigr|^{1/2}.
$$
Again this cannot be $M_0\cap H$-distinguished 
unless $\rho$ is one dimensional. 

\noindent 
(iii) $w=\left(\begin{smallmatrix} 1& & & \\ 
 &0&1& \\ &1&0& \\ & & &1
\end{smallmatrix}\right)$: We have 
$$
M_0\cap wP_0w^{-1}=M_0\cap w^{-1}P_0w=\left\{
\left(\begin{smallmatrix} *&*&0&0 \\ 
0&*&0&0 \\ 0&0&*&* \\ 0&0&0&* 
\end{smallmatrix}\right)\right\}\simeq B_2\times B_2 
$$
where $B_2$ is as in \ref {5.1}. So we may write 
$$
\mathcal F_w([\rho])={\rm Ind}_{B_2\times B_2}^{G_2\times G_2} 
\bigl({}^w([\rho]_{B_2\times B_2})\bigr)  
$$
in this case. 

Now, take $\rho$ to be the Steinberg representation 
${\rm St}_2$ of $G_2$. 
We claim that {\it 
the irreducible $H$-distinguished representation 
$\pi({\rm St}_2)$ is $H$-square integrable}. 
It is enough to look at the exponents coming 
from $M_0\cap H$-distinguished components in 
${\rm I}(\rho)_{P_0}$. So we may look at only the case (iii) above. 
Put $T_2=\left\{ \left(\begin{smallmatrix} *&0\\ 0&*\end{smallmatrix}
\right)\in B_2\right\}$. Under this notation, 
the $T_2\times T_2$-module $[\rho]_{B_2\times B_2}$ is 
given by 
$$
\rho_{B_2}\bigl|\det(\cdot)\bigr|_F^{1/2}
\otimes \rho_{B_2}\bigl|\det(\cdot)\bigr|_F^{-1/2}
=\delta_{B_2}^{1/2}\bigl|\det(\cdot)\bigr|_F^{1/2}
\otimes \delta_{B_2}^{1/2}\bigl|\det(\cdot)\bigr|_F^{-1/2}.  
$$
This is a character of $T_2\times T_2$ written as 
$$
\bigl((t_1,\,t_2),\,(t_3,\,t_4)\bigr)\mapsto \bigl| t_1\bigr|_F\cdot \bigl| t_4\bigr|_F^{-1}. 
$$
Applying $w$ and inducing up to $G_2\times G_2$, 
we have 
\begin{equation} \label{5.2.2} 
\mathcal F_w([{\rm St}_2])={\rm Ind}_{B_2}^{G_2}(|\cdot|_F\otimes 1)\otimes 
{\rm Ind}_{B_2}^{G_2}(1\otimes |\cdot|_F^{-1}). 
\end{equation} 
This is a reducible principal series 
having one dimensional quotient. So the only possible 
element of 
${\mathcal Exp}_{S_0}\bigl(\pi({\rm St}_2)_{P_0}, r_{P_0}(\lambda)\bigr)$ 
(for any $\lambda\in \bigl(\pi(\rho)^*\bigr)^H$) 
is the restriction of the central character 
of \eqref{5.2.2}, which is given by 
$$
\chi\bigl({\rm diag}(s_1,s_1,s_2,s_2)\bigr)=
\bigl|s_1\bigr|_F\cdot \bigl|s_2\bigr|_F^{-1}. 
$$ 
We may conclude that 
$\bigl|\chi(s)\bigr|<1$ for all $s\in S_0^-\setminus Z_0S_0^1$, 
hence the claim follows from Theorem \ref{4.7}.

\end{document}